\documentclass[12pt]{article}
\usepackage{amsmath}
\usepackage{amssymb}
\usepackage{bbm}
\usepackage[margin=2.25cm]{geometry}
\usepackage{hyperref}
\usepackage{tikz}
\usepackage{setspace}

\usetikzlibrary { decorations.pathmorphing, decorations.pathreplacing, decorations.shapes, }

\newtheorem{thm}{Theorem}[section]

\newtheorem{prop}[thm]{Proposition}

\newtheorem{lemma}[thm]{Lemma}

\newtheorem{theorem}[thm]{Theorem}
\newtheorem{remark}[thm]{Remark}
\newtheorem{proposition}[thm]{Proposition}

\newtheorem{corollary}[thm]{Corollary}

\newenvironment{proof}{{\bf Proof:}}{\hfill$\square$\vskip.5cm}
\newenvironment{proofof}{}{\hfill$\square$\vskip.5cm}
\newcommand{\R}{\mathbb{R}}
\newcommand{\N}{\mathbb{N}}

\newcommand{\Z}{\mathbb{Z}

}

\renewcommand{\Z}{{\mathbb{Z}}}

\begin{document} 

\title{Large Deviations for the d'Arcais Numbers}
\author{Shannon Starr}
\date{11 January 2026}

\maketitle

\abstract{
The d'Arcais polynomials $P_n(z)$ for $n\in\{0,1,\dots\}$ are defined as
$\sum_{n=0}^{\infty} P_n(z) q^n = \exp(-z\ln((q;q)_{\infty}))$ where
the $q$-Pochhammer symbol is $(q;q)_{\infty} = \prod_{k=1}^{\infty} (1-q^k)$
for $|q|<1$.
Denoting the coefficients for $n \in \N$ by the formula $P_n(z) = \sum_{k=1}^{n} A(2,n,k) z^k/n!$,
we prove that $k_n! A(2,n,k_n)/n!$ satisfies a Bahadur-Rao type large deviation formula 
in the limit $n \to \infty$ with $k_n/n \to \kappa \in [0,1)$ as long as $k_n \to \infty$.
The large deviation rate function is the Legendre-Fenchel transform
$g^*(-\kappa)$ where $g(\kappa) = f^{-1}(\kappa)$
for the function $f : (0,\infty) \to \R$ given by 
$f(y)= \ln(-\ln((e^{-y};e^{-y})_{\infty}))$.
We relate this fact to information about the abundancy index.}

\section{Introduction}

In a recent {\em tour-de-force}, Abdesselam has derived a local central limit theorem LCLT for the number of commuting $\ell$-tuples of permutations
of $[n]=\{1,\dots,n\}$ with a given size of the orbit when acting on $[n]$, characterized by $k$ \cite{AbdesselamSole}.
More precisely, he considered $n \to \infty$ with $k_n = x\sqrt{n}$ for any $x \in (0,\infty)$.
This is the typicality region because the numbers decay very rapidly beyond this regime.
He has also verified an important conjecture by himself, that the numbers are log-concave in this regime \cite{Abdesselam}.
(Also see \cite{AbdesselamStudents}.)

We are motivated by his work, but we focus our attention on just pairs of commuting permutations, $\ell=2$ because there
are well-known algebraic formulas in this case.
Most notably, the modular symmetry of Dedekind's eta function applies to the $\ell=2$ case.
We believe that Abdesselam may have uncovered new symmetries for the $\ell>2$ case analogous to that.
But we are not yet prepared to make use of his results.

We establish a large deviation formula for the region $k_n/n\to \kappa$ as $n\to\infty$ with $\kappa \in [0,1)$, assuming also $k_n\to\infty$
in case $\kappa=0$.

In case a reader is not familiar with commuting $\ell$-tuples of permutations of $[n]$ with a fixed size of the orbit,
we start the next section with a self-contained definition of the d'Arcais numbers.
In fact, we feel that an equally important connection is the fact that the d'Arcais numbers are the Bell transform of the abundancy index,
which we also introduce in the next section.
We believe that some renewed attention on the abundancy indices could be useful, 
using techniques of analytic number theory.

In Section 6, we will return to discuss the important results of Abdesselam and others, since that is our actual motivation.

\section{Definition of the d'Arcais polynomials}

Consider the following double-series
$$
	1+\sum_{n=1}^{\infty} \sum_{k=0}^{n} \frac{A(2,n,k)}{n!} x^k z^n\,
	=\, \exp\left(-x \ln\big((z;z)_{\infty}\big)\right)\, ,
$$
defined using the $q$-Pochhammer symbol
$$
(q;q)_{\infty}\, =\, \prod_{p=1}^{\infty} (1-q^p)\, ,
$$
defined for $|q|<1$.
%
The polynomials $\mathcal{A}_0(x)\equiv 1$ and 
$$
\forall n \in \N\, ,\
\mathcal{A}_n(x)\, =\, \sum_{k=1}^n A(2,n,k) x^k\, ,
$$
are the d'Arcais polynomials.
The coefficients are the d'Arcais numbers, 
$(A(2,n,k)\, :\, n=\{0,1,\dots\}\, ,\ k \in \{0,\dots,n\})$.
The d'Arcais numbers are the so-called ``Bell transform''
of $n!$ times the abundancy indices
$$
\frac{\sigma(n)}{n}\, =\, \sum_{d|n} \frac{d}{n}\, =\, [z^n]\Big(-\ln\big((z;z)_{\infty}\big)\Big)\, .
$$
Here we use notation from generating functions, as in Wilf \cite{Wilf}.
Being the Bell transform means that
$$
A(2,n,k)\, =\, \frac{n!}{k!}\, [z^n]\bigg(\Big(-\ln\big((z;z)_{\infty}\big)\Big)^k\bigg)\, .
$$
See for example, OEIS for the d'Arcais numbers
$$
\text{\url{https://oeis.org/A008298},}
$$
and the abundancy indices
$$
\text{\url{https://oeis.org/A038048},}
$$
where the Bell transform is stated in the former.

\section{Statement of main results}


Let us define a function $F:(0,\infty) \to (0,\infty)$ as $F(y) = -\ln\big((e^{-y};e^{-y})_{\infty}\big)$.
Note that the first derivative is the Lambert-type series
$$
F'(y)\, =\, -\sum_{n=1}^{\infty} \frac{n e^{-ny}}{1-e^{-ny}}\, ,
$$
and the second derivative is 
$$
F''(y)\, =\, \sum_{n=1}^{\infty} \frac{n e^{-ny}}{(1-e^{-ny})^2}\, .
$$
From all this, it is easy to see that $F$ is surjective. It is strictly decreasing, and it is convex.
Moreover, defining $R_n : (0,\infty) \to (0,\infty)$ for each $n \in \N$ by the formula
$$
R_n(y)\, =\, \frac{\sigma(n)}{n}\, e^{-ny}\, ,
$$
we have
$$
F(y)\, =\, \sum_{n=1}^{\infty} R_n(y)\, ,\qquad
F'(y)\, =\, -\sum_{n=1}^{\infty} n R_n(y)\ \text{ and }\ 
F''(y)\, =\, \sum_{n=1}^{\infty} n^2 R_n(y)\, .
$$
We may define the probability measure $\rho_{\cdot}(y) : \N \to (0,\infty)$ by
$$
\forall n \in \N\, ,\ \text{ we have }\ 
\rho_n(y)\, =\, \frac{R_n(y)}{F(y)}\, .
$$
So then
$$
-\frac{F'(y)}{F(y)}\, =\, \sum_{n=1}^{\infty} n \rho_n(y)\ \text{ and }\
\frac{F''(y)}{F(y)}\, =\, \sum_{n=1}^{\infty} n^2 \rho_n(y)\, .
$$
From this it is easy to see that $-F'(y)/F(y)$ maps $y \in (0,\infty)$ surjectively onto $(1,\infty)$.
Moreover,
$$
\forall y \in (0,\infty)\, ,\ \text{ we have }\
\frac{F''(y)}{F(y)}\, >\, \left(\frac{F'(y)}{F(y)}\right)^2\, .
$$

Our main result is the following local large deviation formula:
\begin{theorem}
\label{thm:BR1}
Suppose $\epsilon>0$ and consider a sequence $k_n$ such that  $k_n/n <1-\epsilon$ for all $n \in \N$,
and  $\lim_{n \to \infty} k_n= \infty$.
Then
$$
\frac{k_n!\, A(2,n,k_n)}{n!}\, \sim\, \frac{e^{k_n \ln(F(y_n)) + n y_n}}
{\sqrt{2\pi  k_n \left(\frac{F''(y_n)}{F(y_n)} - \left(\frac{F'(y_n)}{F(y_n)}\right)\right)}}\, ,\
\text{ as $n \to \infty$,}
$$
where for each $n$ we define $\kappa_n=k_n/n$ and we let $y_n \in (0,\infty)$ be the unique solution of 
$$
-\frac{F'(y_n)}{F(y_n)}\, =\, \frac{1}{\kappa_n}\, .
$$
\end{theorem}
\begin{figure}
\begin{center}
\begin{tikzpicture}
\draw (0,0) node[] {\includegraphics[width=11cm]{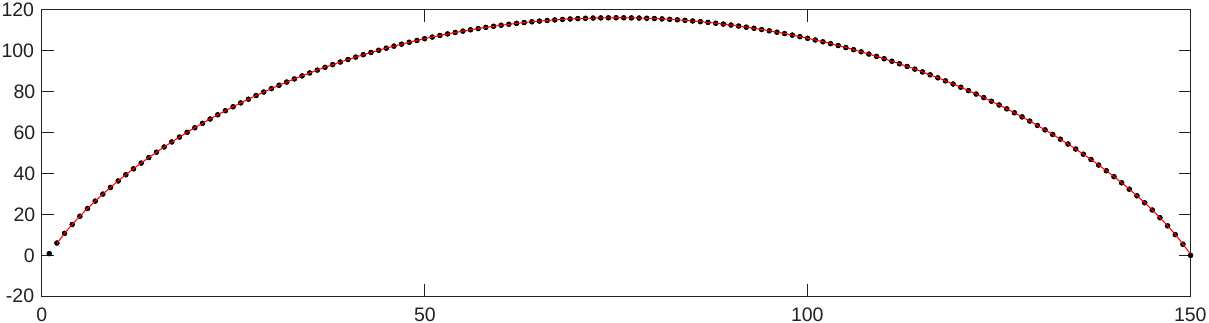}};
\draw (9,0) node[] {\includegraphics[width=5cm]{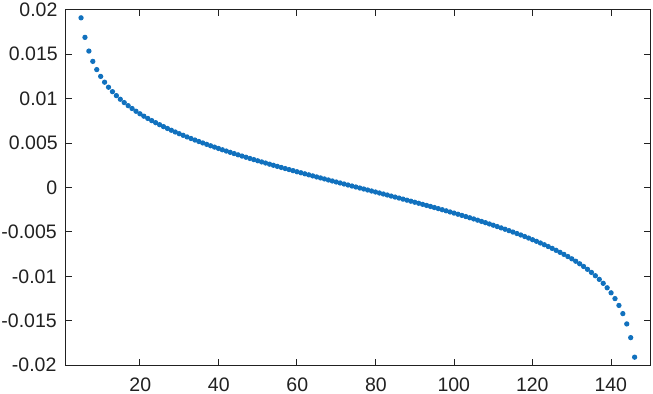}};
\end{tikzpicture}
\caption{
In the first plot, the blue dots are a list-plot of $\ln(a(n,k)) = \ln(k! A(2,n,k)/n!)$ versus $k$, for $n=150$ and $k \in \{1,\dots,150\}$.
The red curve is the logarithm of the approximation from Theorem \ref{thm:BR1}.
In the second plot, in blue, we have plotted $(a(n,k)-a(n,n+1-k))/(a(n,k)+a(n,n+1-k))$ for $n=150$ and for $k$ from $3$ to $147$.
\label{fig:antialias}
}
\end{center}
\end{figure}
In Figure \ref{fig:antialias}, in the first plot, we demonstrate the basic agreement between the numbers and their 
approximation.
The plot is {\em not} perfectly symmetric about the midpoint, as we will demonstrate in the next section.
Let us denote
$$
a(n,k)\, =\, \frac{k! A(2,n,k)}{n!}\, .
$$
Then we may measure the asymmetry as $(a(n,k)-a(n,n+1-k))/(a(n,k)+a(n,n+1-k))$.
We will prove that this does not vanish in the large $n$ limit.

Our motivation for all of these results is Abdesselam's LCLT, which we will discuss more in Section 6.
Abdesselam used his LCLT to prove log-concavity for large $n$, in the typicality region, where $k_n = x \sqrt{n}$ for $x \in (0,\infty)$.
That is the region where the large deviation rate function is largest.
But we establish the same log-concavity result for the region where our asymptotic formula holds true.
\begin{corollary}
\label{cor:logConcave}
Let us define $a(n,k) = k! A(2,n,k)/n!$. Let us define $\mathcal{K}(y) = -F(y)/F'(y)$ so that $\mathcal{K}(y_n) = \kappa_n$.
Then with the same set-up as above, we have
\begin{equation}
\label{eq:logConcave}
\ln\left(\frac{\big(a(n,k_n)\big)^2}{a(n,k_n-1)a(n,k_n+1)}\right)\, \sim\, \frac{1}{n \big(\mathcal{K}(y_n)\big)^3 \mathcal{V}(y_n)}\, ,
\end{equation}
as $n \to \infty$, where
$$
\mathcal{V}(y)\, =\, \frac{F''(y)}{F(y)} - \left(\frac{F'(y)}{F(y)}\right)^2\, .
$$
In particular, the ratio is positive: the sequence $a(n,k)$ is log-concave for sufficiently large $n$, assuming $k_n/n \in [0,1-\epsilon]$
and assuming $k_n \to \infty$.
\end{corollary}
\begin{figure}
\begin{center}
\begin{tikzpicture}
\draw (0,0) node[] {\includegraphics[width=14cm]{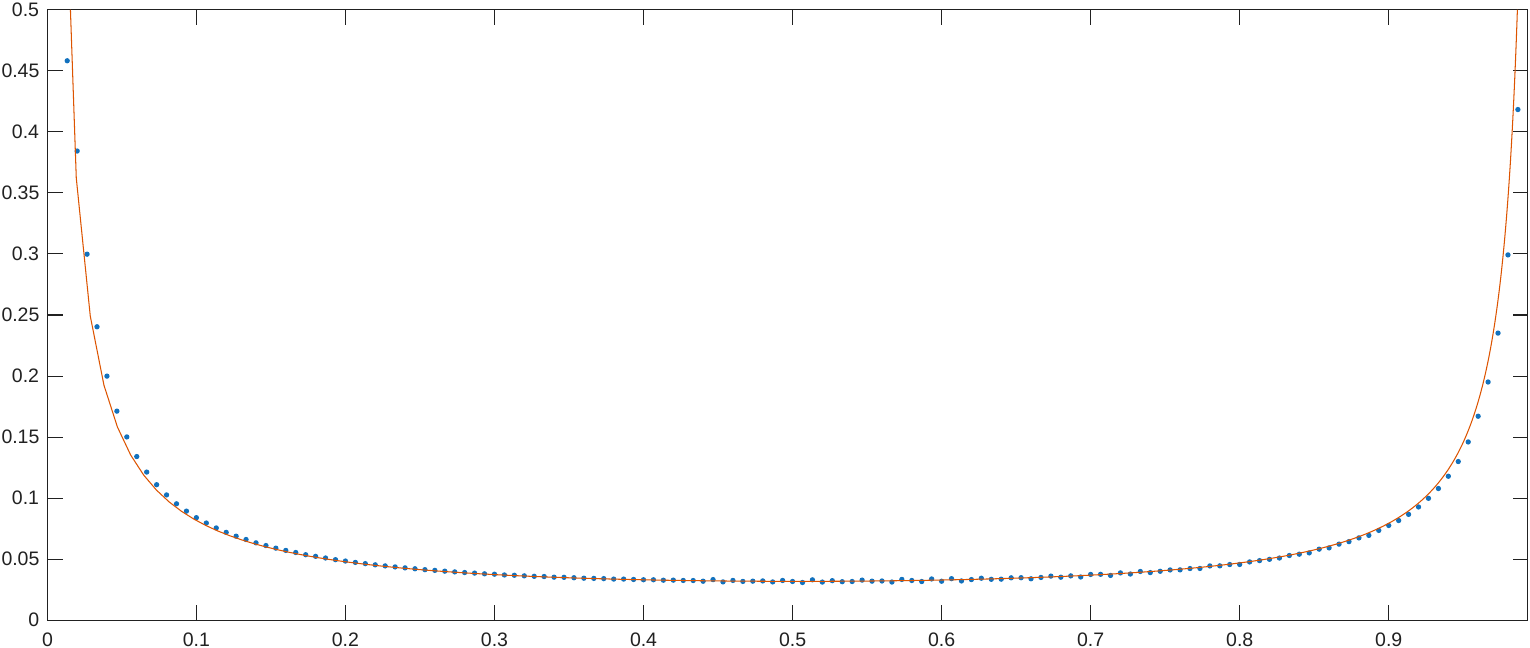}};
\end{tikzpicture}
\caption{
In the  plot, the blue dots are $\ln\Big(\big(a(n,k)\big)^2/\big(a(n,k-1)a(n,k+1)\big)\Big)$ versus $k/n$, for $n=150$ and $k \in \{2,\dots,149\}$.
The red curve is a plot of $1/\Big(n \big(\mathcal{K}(y)\big)^3 \mathcal{V}(y)\Big)$ versus $\mathcal{K}(y)$ for a range of $y$ between $0.01$ and $10$.
\label{fig:logConcave}
}
\end{center}
\end{figure}
In Figure \ref{fig:logConcave}, we compare the logarithm of the ratio on the left-hand-side of (\ref{eq:logConcave})
to the curve $\bigg(\mathcal{K}(y)\, ,\ 1/\Big(n \big(\mathcal{K}(y)\big)^3 \mathcal{V}(y)\Big)\bigg)$ from the right-hand-side for $n=150$.

\section{Brief observations about the main result}
We wish to make two brief observations about Theorem \ref{thm:BR1}.
Firstly, the large deviation rate function is not symmetric (although the asymmetry is somewhat small).
Secondly, the formula is consistent with the Bahadur-Rao asymptotic formula for local large deviations.

\subsection{Proof of asymmetry of the large deviation rate function}
\begin{lemma}
\label{lem:asym}
We have
$$
\lim_{y \to \infty} F(y) F(y+\ln(F(y)))\, =\, \frac{2\zeta(2)}{3}\, =\, \frac{\pi^2}{9}\, \approx\, 1.09662\, .
$$
Moreover, while $\lim_{y \to 0^+} F(y) F(y+\ln(F(y)))=1$, we have
$$
F(y) F(y+\ln(F(y)))-1\, \sim\, -\frac{\pi^2-9}{\pi^2}\, y\, , \text{ as  $y \to 0^+$,}
$$
\end{lemma}
\begin{proofof}{\bf Proof of Lemma \ref{lem:asym}:}
Note that
using Dedekind's eta function, 
$$
\eta(\tau)\, =\, q^{1/24} (q;q)_{\infty}\bigg|_{q=\exp\left(2\pi i \tau\right)}\, =\, e^{\pi i \tau/12} (e^{2\pi i \tau};e^{2\pi i \tau})_{\infty}\, ,
$$
we have
$$
F(y)\, =\, -\ln(\eta(iy/2\pi))-\frac{y}{24}\, .
$$
This is particularly useful for asymptotics when $q \to 1$ due to the modular functional relation:
$$
\eta\left(-\frac{1}{\tau}\right)\, =\, \sqrt{-i\tau}\, \eta(\tau)\, ,
$$
where $\sqrt{-i\tau}=1$ for $\tau=i$.
This has been noted before, for example in Moak \cite{Moak}. 

Using the modular symmetry,
$$
F(y)\, =\, \frac{\zeta(2)}{y} - \frac{\ln(y)}{2} + \frac{\ln(2\pi)}{2} - \frac{y}{24} + F\left(\frac{4\pi^2}{y}\right)\, .
$$
We notice that using this formula for $0<y\ll 1$, 
one can think of $F(4\pi^2/y)$ as a correction term.
But as such, it is asymptotically smaller than any power of $y$ as $y \to \infty$.
This leads to asymptotic formulae which appear to be exact after a finite number of terms (even though the actual functions have these essentially singular corrections).
A good example is in Banerjee and Wilkerson \cite{BanerjeeWilkerson} whose formula for $\ln(q^x;q)_{\infty}$, when specialized to $x=1$, terminates after the second term
because of the appearance of a factor in the coefficients $\zeta(-n)B_n$ where $\zeta(-n)=0$ for $n$ even and  the Bernoulli numbers $B_n=0$ for $n$ odd.
But Banerjee and Wilkerson explicate a complete asymptotic series, with non-vanishing terms, when $x$ is not 1.

For large $y$, we have $F(y) = e^{-y} + \frac{3}{2} e^{-2y} + O(e^{-3y})$
and $\ln(F(y)) = -y + \frac{3}{2}\, e^{-y} + O(e^{-2y})$.
So $y+\ln(F(y)) = \frac{3}{2}\, e^{-y} + O(e^{-2y})$. Taking $v = \frac{3}{2}\, e^{-y}$, which is small, the modular symmetry
gives $F(v) \sim 2 \zeta(2)  e^{y}/3$. Thus $F(y)F(v) \to 2 \zeta(2)/3$.

For $y$ small, $y+\ln(F(y))$ is large: it is $y + \ln\left(\frac{\zeta(2)}{y} - \frac{\ln(y)}{2}+\frac{\ln(2\pi)}{2}-\frac{y}{24}\right)$
to all orders in a complete asymptotic series. Defining this to be $v$, one can obtain an asymptotic expansion of $F(v)$.
Then to all orders in a complete asymptotic expandion we have $\left(\frac{\zeta(2)}{y} - \frac{\ln(y)}{2}+\frac{\ln(2\pi)}{2}-\frac{y}{24}\right) F(v)$
for $F(y) F(y+\ln(F(y)))$.
Then a standard Taylor expansion gives the second result. 
\end{proofof}
The reason that this lemma implies asymmetry of the LDP rate function is via the contrapositive.
Symmetry of the LDP rate function would imply the result: $F(y) F(y+\ln(F(y)))$ would be equal to $1$, identically.
We will demonstrate this now.

Define $\mathcal{Y} : (0,1) \to (0,\infty)$ such that for each $\kappa$, we let $\mathcal{Y}(\kappa)$ be the unique solution $y$ of 
$$
-\frac{F'(y)}{F(y)}\, =\, 1/\kappa\, .
$$
Then the large deviation rate function (modulo a sign) is $\Gamma : (0,1) \to (0,\infty)$ defined by
$$
\Gamma(\kappa)\, =\, \kappa \ln(F(\mathcal{Y}(\kappa))) + \mathcal{Y}(\kappa)\, .
$$
This is also the same as 
$$
\Gamma(\kappa)\, =\, \min_{y \in (0,\infty)} \left(\kappa \ln(F(y))+y\right)\, ,
$$
which can be related to the Legendre transform.
Define the strictly convex function $f : (0,\infty) \to \R$ by
$$
f(y)\, =\, \ln(F(y))\, .
$$
(Note that $f''(y) = (F''(y)/F(y)) - (F'(y)/F(y))^2$ which we already noted is positive.)
This is strictly decreasing and surjective onto $\R$.
Then 
$$
\Gamma(\kappa)\, =\, \min_{y \in (0,\infty)} \left(y + \kappa f(y)\right)\, = -\kappa \max_{y \in (0,\infty)} \left(-\frac{1}{\kappa}\, y - f(y)\right)\, .
$$
So, using the Legendre transform, we have $\Gamma(\kappa) = -\kappa f^*(-1/\kappa)$.
The idea is to use the fact that if the Legendre transform of a strictly convex function were even,
then by duality of the Legendre transform, that would imply that the original function was even.

To use this more effectively, let us define $g : \R \to (0,\infty)$ to be the inverse-function of $f$.
Then we have
$$
\Gamma(\kappa)\, =\, \min_{u \in \R} \left(g(u) + \kappa u\right)\, =\, -\max_{u \in \R}\left(-\kappa u-g(u)\right)\, .
$$
So, in terms of the Legendre transform, $g^*$, we have
$$
\Gamma(\kappa)\, =\, -g^*(-\kappa)\, .
$$
Now define $\widetilde{\Gamma}(\lambda) = \Gamma(\frac{1}{2}+\lambda)$.
The claim is that $\widetilde{\Gamma}$ is not an even function.
If it were, we would have $\widetilde{\Gamma}(\lambda) =-g^*(-\frac{1}{2}-\lambda) = \widetilde{\Gamma}(-\lambda) = -g^*(-\frac{1}{2}+\lambda)$.
Note that
$$
g^*\left(-\frac{1}{2}+\lambda\right)\, =\, \max_{u \in \R} \left(\lambda u - \frac{1}{2}\, u -g(u)\right)\, .
$$
So, defining $\widetilde{g}(u)\, =\, \frac{1}{2}\, u + g(u)$, we have
$$
\widetilde{\Gamma}(\lambda)\, =\, -\widetilde{g}^*(-\lambda)\, .
$$
Then if $\widetilde{\Gamma}$ were even that would imply $\widetilde{g}^*$ is even.
But that would imply $\widetilde{g}$ is even because $\widetilde{g}$ is the Legendre transform of $\widetilde{g}^*$
by Legendre duality for strictly convex functions.
(And the Legendre transform of an even function is itself an even function.)
In turn that would mean
$$
g(-u)\, =\, g(u) + u\, .
$$
Taking $u=f(y)$ that would mean
$$
g(-f(y))\, =\, y+f(y)\, .
$$
Then applying $f$ to both sides
$$
-f(y)\, =\, f(y+f(y))\, .
$$
Since $f(y)=\ln(F(y))$, that means $1/F(y)=F(y+\ln(F(y)))$. 
So symmetry would imply that $F(y)F(y+\ln(F(y)))$ would be identically 1.

\begin{figure}
\begin{center}
\begin{tikzpicture}
\draw (0,0) node[] {\includegraphics[width=15cm]{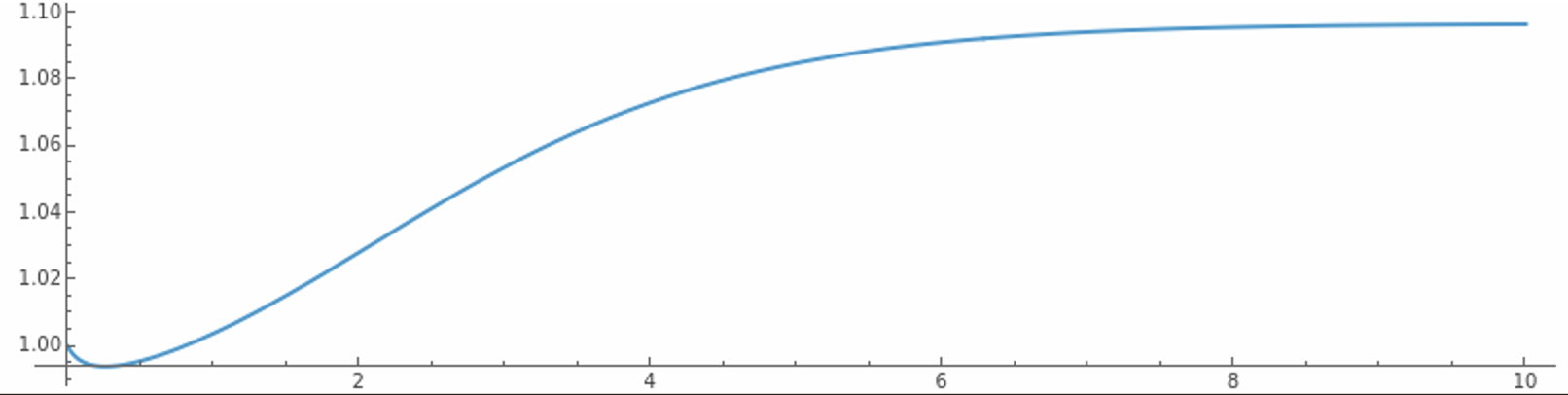}};
\end{tikzpicture}
\caption{
Here we plot $F(y) F(y+\ln(F(y)))$ versus $y$. 
\label{fig:asym}
}
\end{center}
\end{figure}

But by the lemma, we know it is not identically 1.
We plot $F(y)F(y+\ln(F(y)))$ in Figure \ref{fig:asym}.
It is somewhat close to 1. But it is not identically 1.

\subsection{Comparison to the Bahadur-Rao theorem}

Suppose $\mathsf{X}_1,\mathsf{X}_2,\dots$ are independent, identically distributed
integer-valued random variables where the distribution has step-size 1.
(So, defining $\operatorname{supp}(\mathsf{X}_1)$
to be the set  
$\{x \in \Z\, :\, \mathbf{P}(\mathsf{X}_1=x)>0\}$,
we have $\operatorname{gcd}(\{x_1-x_2\, :\, x_1-x_2 \in \operatorname{supp}(\mathsf{X}_1)\})$
equals 1.)
Let
$$
\mathcal{L}(\lambda)\, =\, \ln\left(\mathbf{E}\left[e^{\lambda \mathsf{X}_1}\right]\right)\, ,
$$
be the logarithmic moment generating function (cumulant generating function)
which is strictly convex.
\begin{proposition}
For an integer sequence $k_n$ such that $\kappa_n = k_n/n$ remains bounded,
$$
\mathbf{P}\left(\mathsf{X}_1+\dots+\mathsf{X}_n=k_n\right)\,
\sim\, \frac{\exp\left(-n \mathcal{L}^*(\kappa_n)\right)}{\sqrt{2\pi n \sigma_{n}^2}}\, ,
$$
where the Legendre transform is $\mathcal{L}^*(\kappa_n) = s_n \kappa_n - \Lambda(s_n)$,
where $s_{n}$ is the choice of $s$ that maximizes
$$
s \kappa_n  - \mathcal{L}(s)\, ,
$$
(so that $\mathcal{L}'(s_n)=\kappa_n$)
and  
$$
\sigma_n^2\, =\, \mathcal{L}''(s_n)\, .
$$
\end{proposition}
See, for example, Theorem 3.7.4 in Dembo and Zeitouni \cite{DemboZeitouni},
although the version stated here is a local version, with the lattice setting as above.

Suppose in the setting of the Bahadur-Rao theorem, we define the generating function
$$
\mathcal{G}_{\mathsf{X}_1}(z)\, =\, \sum_{n=-\infty}^{\infty} z^n \mathbf{P}(\mathsf{X}_1=n)\, .
$$
Then
$$
\mathbf{P}\left(\mathsf{X}_1+\dots+\mathsf{X}_n=k\right)\,
=\, [z^k] \Big(\mathcal{G}_{\mathsf{X}_1}(z)^n\Big)\, .
$$
Compare this to the formula that we use for the proof of Theorem \ref{thm:BR1}:
$$
\frac{k! A(2,n,k)}{n!}\, =\, [z^n]\bigg(\Big(-\ln\big((z;z)_{\infty}\big)\Big)^k\bigg)\, .
$$
The positions of $n$ and $k$ are switched relative to the Bahadur-Rao theorem.

In order to make the result  of Theorem \ref{thm:BR1} appear closer to the Bahadur-Rao theorem, let us define
a function $\gamma : \R \to (0,1)$ as 
$$
\gamma(u)\, =\, g(-u)\, .
$$
This is still strictly convex. But now $g^*(-\kappa)=\gamma^*(\kappa)$.
Right before Theorem \ref{thm:BR1}, we noted that $-F'/F$ maps surjectively onto $(1,\infty)$.
From this we know that $f=\ln(F)$ has derivative mapping surjectively to $(-\infty,-1)$.
So, by the inverse mapping formula, $g$ has derivative mapping surjectively onto $(-1,0)$.
So $\gamma$ has derivative mapping surjectively onto $(0,1)$.
And since $\gamma$ is strictly convex, its derivative is strictly increasing on $\R$.
So Theorem \ref{thm:BR1} may be rewritten as follows:
as long as $k_n/n \in [0,1-\epsilon)$ for all $n$ and $k_n \to \infty$, then we have
\begin{equation}
\label{eq:BR2}
\frac{k! A(2,n,k)}{n!}\, \sim\, \frac{e^{-n\gamma^*(\kappa_n)}}{\sqrt{2\pi n \sigma^2(\kappa_n)}}\, \cdot \kappa_n\, ,
\end{equation}
as $n \to \infty$, where $\kappa_n = k_n/n$ and for each $\kappa \in (0,1)$, if we let $\mathcal{U}(\kappa)$ be the unique $u$ such that $\gamma'(u)=\kappa$
then $\sigma^2(\kappa) = \gamma''(\mathcal{U}(\kappa))$.

Relative to the expected formula from the Bahadur-Rao formula, in our formula there appears an extra factor of 
$$
\kappa_n\, =\, \gamma'(\mathcal{U}(\kappa_n))\, .
$$
The reason for this is that we are ultimately using a contour integral
via Cauchy's integral formula
\begin{equation}
\label{eq:Cauchy}
\frac{k! A(2,n,k)}{n!}\, =\, 
\int_{\mathcal{C}(0;\rho)} z^{-n}\exp\left(k\ln\left(-\sum_{p=1}^{\infty}\ln(1-z^p)\right)\right)\, \frac{dz}{2\pi i z}\, ,
\end{equation}
for $\rho=e^{-y^*}$ for $y^*$ solving $f'(y^*)=-1/\kappa$.
Then $z=e^{-(y^*-i\theta)}$,  for $-\pi<\theta<\pi$.
We can transform the variables to more closely resemble the Bahadur-Rao formula
by taking $y-i\theta = \gamma(u-i\phi)$.
But then $dz/z = id\theta = i\gamma'(u-i\phi) d\phi$.
So, taking $u=\mathcal{U}(\kappa)$ and expanding around $\phi=0$ we have an extra factor of $\gamma'(\mathcal{U}(\kappa))$.

\section{Proof of Theorem \ref{thm:BR1}}

The proof is elementary, based on the Cauchy integral formula.
More specifically, we use the Hayman version of the saddle point method.
An excellent reference is the textbook of Flajolet and Sedgewick.
But we will try to be self-contained, here.

Consider the Cauchy integral formula
$$
a(n,k)\, =\, e^{ny+kF(y)} \int_{-\pi}^{\pi} e^{-in\theta} \exp\Big(k\ln\big(F(y-i\theta)\big)-k\ln\big(F(y)\big)\Big)\, \frac{d\theta}{2\pi}\, .
$$
The idea is to Taylor expand the integrand (as a function of $\theta$) in the region of $\theta \in (-\tau_n,\tau_n)$ where 
$\tau_n>0$ is determined by choosing
$$
k\ln\big(F(y-i\theta)\big)-k\ln\big(F(y)\big)-in\theta
$$
to be on the order of 1, there.

But equally important, one must obtain {\em uniform} bounds that apply also for $\theta \in (-\pi,\pi) \setminus (-\tau_n,\tau_n)$.
The uniform bounds are supposed to be such that one may determine that 
$$
\int_{(-\pi,\pi) \setminus (-\tau_n,\tau_n)} e^{-in\theta} \exp(\cdots)\, \frac{d\theta}{2\pi}\,
\ll\,
\int_{(-\tau_n,\tau_n)} e^{-in\theta} \exp(\cdots)\, \frac{d\theta}{2\pi}\, ,
$$
as $n \to \infty$.
(And of course the goal in the $(-\tau_n,\tau_n)$ regime is to establish that $\int_{(-\tau_n,\tau_n)} (\cdots)\, d\theta$ equals
the stated asymptotic formula of the theorem.)

In this context, following the original circle method of Hardy and Ramanujan, the regime $(-\tau_n,\tau_n)$ is called the ``major arc,'' and $(-\pi,\pi) \setminus (-\tau_n,\tau_n)$ is called the ``minor arc.''
We are thinking of these as subsets of the circle $\R / (2\pi \Z)$.

The main issue is to control the bounds when $y>0$ is small.
So let us state two lemmas to demonstrate how this is controlled.

\subsection{Two lemmas}
\begin{lemma}
\label{lem:Taylor}
For $y \in (0,\infty)$ and 
$\vartheta \in \R$,
\begin{equation*}
\begin{split}
\frac{iy}{\mathcal{K}(y)}\, \vartheta + \ln(F(y-iy\vartheta))\, - \ln(F(y))\, 
=\, 
\sum_{r=2}^{R}\frac{(-i)^r}{r!}\, {\mathcal{W}}^{(r)}(y)  \vartheta^r + \mathcal{E}_R(y;\vartheta)\, ,
\end{split}
\end{equation*}
where 
$$
\mathcal{W}^{(r)}(y)\, =\,  y^rf^{(r)}(y)\, ,
$$
for $r=3,\dots,R$, and
$$
\mathcal{E}_R(y;\vartheta)\, =\, \frac{(-i)^{r+1} y^{r+1} \vartheta^{r+1}}{r!} \int_0^1 (1-s)^r f^{(r+1)}(y-iy\vartheta s)\, ds\, .
$$
These are all analytic functions of $y \in \{z\, :\, \operatorname{Re}(z)>0\}$ and $\vartheta \in \{z\, :\, \operatorname{Im}(z)>-1\}$.
Moreover, we have
$\lim_{y \to 0^+} \mathcal{W}^{(r)}(y) = (-1)^r(r)!$, while
$$
\lim_{y \to 0^+} \mathcal{E}_R(y;\vartheta)\, =\, -\ln(1-i\vartheta)-\sum_{r=1}^{R} \frac{i^r \vartheta^r}{r}\, , 
$$
which is clearly $O(|\vartheta|^{R+1})$.
\end{lemma}

The point of the first lemma is to consider the Taylor expansion at this scale, and establish uniformity in $y$.
We establish that the coefficients  are uniformly bounded if we assume $y \in (0,B]$ for some $B<\infty$.
Moreover, the error term is such that there exists a uniformly bounded constant $C<\infty$ (which is a function $C(B)$)
and a $\vartheta_0>0$ (also depending on $B$)
such that $|\mathcal{E}_R(y;\vartheta)|\leq C|\vartheta|^{R+1}$ for $|\vartheta|\leq \vartheta_0$.

To treat the regime outside the small interval, we need separate uniform bounds.

\begin{lemma}
\label{lem:Uniform}
We have, for all $y>0$ and $\theta \in \R$,
$$
1 - \frac{|F(y-i\theta)|^2}{\big(F(y)\big)^2}\, \geq\, \beta(y,\theta)\, \stackrel{\mathrm{def}}{:=}\, 
\frac{F(2y)}{\big(F(y)\big)^2}\, \cdot \frac{\cosh(y/2) \sin^2(\theta/2)}{\sinh(y/2)\big(\sinh^2(y/2)+\sin^2(\theta/2)\big)}\, .
$$
Moreover, rescaling we have
$$
\beta(y,y\vartheta)\, \geq\, \frac{2F(2y)}{y\big(F(y)\big)^2 } \cdot \frac{\vartheta^2}{\pi^2 y^{-2} \sinh^2(y/2)+\vartheta^2}\, .
$$
\end{lemma}
Let us recall that $F(y) \sim \zeta(2)/y$ as $y \to 0^+$. Thus we have $\beta(y,y\vartheta) \sim 4 \vartheta^2 / (1+\vartheta^2)$.
The second lemma is such that 
$$
\left|\exp\left(-k\, \ln\left(\frac{F(y-i\theta)}{F(y)}\right)\right)\right|\, \leq\, \exp\left(-k\beta(y,\theta)/2\right)\, ,
$$
because $-\ln(1-x)\geq x$.
So if $\theta = y \vartheta$ for $|\vartheta| > \vartheta_0$,
then we have a bound 
$$
\left|\int_{(-\pi,\pi) \setminus (-y\vartheta_0,y\vartheta_0)} \left(\cdots\right)\, \frac{d\theta}{2\pi}\right|\, 
\leq\, \exp\left(-k\beta(y,y\vartheta_0)\right)\, .
$$
We note that $\beta(y,\theta)$ is an even function of $\theta$, and such that it is strictly increasing on $\theta \in (0,\pi)$.
Since we assume $k$ is a large parameter, this type of bound suffices for our purpose.
\begin{figure}
\begin{center}
\begin{tikzpicture}
\draw (0,0) node[] {\includegraphics[width=12cm]{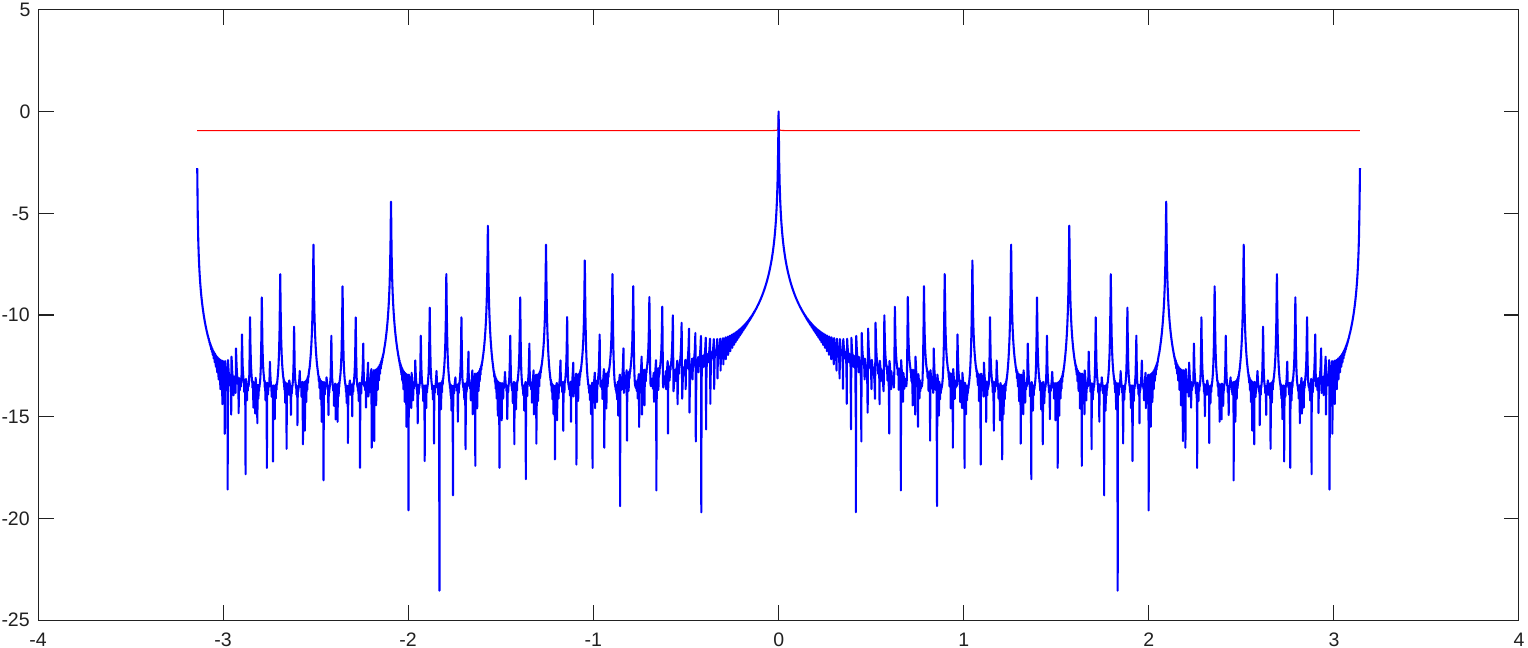}};
\draw (7.5,0) node[] {\includegraphics[width=5cm]{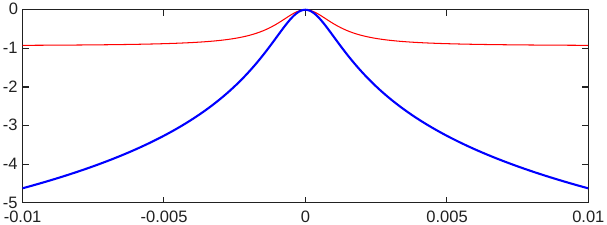}};
\end{tikzpicture}
\caption{
Here we plot $\ln\Big(\big|F(y+i\theta)\big|^2/\big(F(y)\big)^2\Big)$ for $y=10^{-3}$. We truncated the sum defining $F$ at $n=10^5$.
In red is the logarithm of the bound $1-\beta(y;\theta)$. The inset is a zoom-in near the maximum peak, at $\theta=0$. 
\label{fig:Farey}
}
\end{center}
\end{figure}
In Figure \ref{fig:Farey} we plot $\ln\Big(\big|F(y+i\theta)\big|^2/\big(F(y)\big)^2\Big)$ versus $\theta$ for a choice of $y=10^{-3}$.
This function does not decay monotonically.
But $\beta(y,\theta)$ does result in a monotonically decaying uniform bound.

\subsection{Proofs of the lemmas}

\begin{proofof}{\bf Proof of Lemma \ref{lem:Taylor}:}
The first part of the lemma is the statement of Taylor's theorem with remainder, proved by the fundamental theorem of calculus and induction.
For the second part of the lemma, note that 
$$
F(y)\, =\, -\sum_{n=1}^{\infty} \ln\left(1-e^{-ny}\right)\, .
$$
Thus we have
$$
yF(y)\, =\, \int_0^{\infty} \sum_{n=1}^{\infty} \mathbf{1}_{[(n-1)y,ny)}(t) (-\ln(1-e^{-ny}))\, dy\, ,
$$
where the integrand is dominated by $-\ln(1-e^{-t})$ by monotonicty, and the latter integrand is positive and integrable.
So by DCT
$$
\lim_{y \to 0^+} y F(y)\, =\, \int_0^{\infty} (-\ln(1-e^{-t}))\, dt\, =\, \zeta(2)\, .
$$
Hence $f(y) \sim \ln(y)$ as $y \to 0^+$.
Now
$$
F^{(r)}(y)\, =  \sum_{n=1}^{\infty} n^r \frac{d^r}{dz^r}\Big(-\ln(1-e^{-z})\Big)\, |_{z=yn}\, .
$$
So, by a Riemann sum approximation of a Riemann integral, we have
$$
\lim_{y \to 0^+} y^{r+1} F^{(r)}(y)\, =\, \int_0^{\infty} z^r\, \frac{d^r}{dz^r}\Big(-\ln\left(1-e^{-z}\right)\Big)\, dz\, .
$$
But by integration-by-parts applied $r$ times, this is $(-1)^r \zeta(2) (r!)$.
Then one may use the relation between cumulants and moments to determine the desideratum.
This is also known as Fa\`a de Bruno's formula, and uses the incomplete Bell polynomials.
\end{proofof}

Instead of the argument above given for the second part of the lemma, where one calculates the leading order behavior of the Taylor coefficients
as $y \to 0^+$, one can instead use the analyticity of $F$ and the modular symmetry, resulting in 
$$
F(y)\, =\, \frac{\zeta(2)}{y} - \frac{\ln(y)}{2} + \frac{\ln(2\pi)}{2} - \frac{y}{24} + F\left(\frac{4\pi^2}{y}\right)\, .
$$
This is more direct. We gave the indirect argument above, in hopes that it might be useful at a later date
to generalize the result beyond the LDP for d'Arcais numbers to the more general numbers $A(\ell,n,k)$ for $\ell>2$.

\begin{proofof}{\bf Proof of Lemma \ref{lem:Uniform}:}
Recall that one may write
$$
F(y)\, =\, \sum_{n=1}^{\infty} \frac{\sigma(n)}{n}\, e^{-ny}\, .
$$
We have previously define
$$
R_n(y)\, =\, \frac{\sigma(n)}{n}\, e^{-ny}\ \text{ and }\ 
\rho_n(y)\, =\, \frac{R_n(y)}{F(y)}\, .
$$
So
$$
\frac{F(y-i\theta)}{F(y)}\, =\, \sum_{n=1}^{\infty} \rho_n(y) e^{in\theta}\, ,
$$
which is the Fourier transform of the probability mass function $n\mapsto\rho_{n}(y)$.
In particular,
$$
\frac{|F(y-i\theta)|^2}{\big(F(y)\big)^2}\, =\, \sum_{m=1}^{\infty} \sum_{n=1}^{\infty} \rho_m(y) \rho_n(y) e^{i(n-m)\theta}\, 
=\, \sum_{m=1}^{\infty} \sum_{n=1}^{\infty} \rho_m(y) \rho_n(y) \cos\big((n-m)\theta\big)\, ,
$$
where we used symmetry at the last step. Then, using the half-angle formula $\cos(\phi) = 1 - 2\sin^2(\phi/2)$, and the fact that $\sum_{n=1}^{\infty} \rho_n(y)=1$,
we have
$$
\frac{|F(y-i\theta)|^2}{\big(F(y)\big)^2}\, =\, 1 - 2\sum_{m=1}^{\infty} \sum_{n=1}^{\infty} \rho_m(y) \rho_n(y) \sin^2\left(\frac{(n-m)\theta}{2}\right)\, .
$$
In turn, this equals
$$
\frac{|F(y-i\theta)|^2}{\big(F(y)\big)^2}\, =\, 1 - 4\sum_{n=1}^{\infty} \sum_{k=1}^{\infty} \rho_n(y) \rho_{n+k}(y) \sin^2\left(\frac{k\theta}{2}\right)\, .
$$
But now, rearranging
$$
1-\frac{|F(y-i\theta)|^2}{\big(F(y)\big)^2}\, =\, 4\sum_{n=1}^{\infty} \sum_{k=1}^{\infty} \rho_n(y) \rho_{n+k}(y) \sin^2\left(\frac{k\theta}{2}\right)\, ,
$$
all the summands in the final summation are nonnegative.
So we may use monotonicity.
We note that $\rho_{n+k}(y) = \big(F(y)\big)^{-1}e^{-ny} e^{-ky} \sigma(n+k)/(n+k)$.
But $\sigma(n+k)/(n+k) \geq 1$.
So we have
$$
1-\frac{|F(y-i\theta)|^2}{\big(F(y)\big)^2}\, \geq\, \frac{4}{F(y)}\, \sum_{n=1}^{\infty} \sum_{k=1}^{\infty} e^{-ny}\rho_n(y) e^{-ky} \sin^2\left(\frac{k\theta}{2}\right)\, .
$$
Splitting the two summations, now, gives the result.
\end{proofof}

\subsection{Completion of the proof}

In the Taylor expansion, in Lemma \ref{lem:Taylor}, if $y \in (0,B]$ then rescaling $\theta=y\vartheta$
we have established sufficient uniformity  to guarantee that there exists a $\vartheta_0>0$ such that 
$$
\frac{iy}{\mathcal{K}(y)}\, \vartheta + \ln(F(y-iy\vartheta))\, - \ln(F(y))\, =\, -\frac{1}{2}\, \mathcal{W}^{(2)}(y) \vartheta^2 + O(\vartheta^3)\, ,
$$
and in particular
$$
\operatorname{Re}[\ln(F(y-iy\vartheta))]\, - \ln(F(y))\, \leq\, -\frac{1}{4}\, \mathcal{W}^{(2)}(y) \vartheta^2\, ,
$$
for all $\vartheta \in [-\vartheta_0,\vartheta_0]$.

We are not done with the integral over $\theta \in (-y\vartheta_0,y\vartheta_0)$ because we have to further rescale, due to the large parameter $k$,
in order to focus on an even smaller interval.
But before doing that, let us note that the complementary integral is negligible.

Due to Lemma \ref{lem:Uniform}, we have
$$
\left|\int_{(-\pi,\pi)\setminus (-\vartheta_0,\vartheta_0)} e^{-in\theta} \exp\left(k\ln\big(F(y-i\theta)\big)-k\ln\big(F(y)\big)\right)\, \frac{d\theta}{2\pi}\right|\,
\leq\, e^{(k/2)\ln(1-\beta(y,y\vartheta_0))}\, ,
$$
and this may be further bounded above by $e^{-k \beta(y,y\vartheta_0)/2}$.
But since $y \in (0,B]$, there is an absolute constant $b>0$ such that $\beta(y,y\vartheta_0)\geq b$ throughout the interval.
So that integral is bounded above by $e^{-bk/2}$ which vanishes faster as $k \to \infty$ than any algebraic power of $1/k$.

Now returning to the integral in the ``major arc'' $\theta \in (-y\vartheta_0,y\vartheta_0)$, we have that the integrand
is asymptotically $\exp(-k\mathcal{W}^{(2)}(y)\vartheta^2/2)$.
Moreover, the magnitude of the integrand is bounded by $\exp(-k\mathcal{W}^{(2)}(y)\vartheta^2/4)$.
Let us rescale $\vartheta = k^{-1/2} \Theta$.
Then we obtain an integral
$$
\int_{-k^{1/2} \vartheta_0}^{k^{1/2} \vartheta_0} \mathrm{Integrand}(\vartheta)\Bigg|_{\vartheta = k^{-1/2} \Theta}\, \frac{yk^{-1/2} d\Theta}{2\pi}\, .
$$
But now the integrand after the substitution becomes asymptotically $\exp(-\mathcal{W}^{(2)}(y)\Theta^2/2)$
and bounded in magnitude by $\exp(-\mathcal{W}^{(2)}(y)\Theta^2/4)$.
So the result follows by the dominated convergence theorem.
{\bf Q.E.D.}

\subsection{Considerations of the abundancy index asymptotics}

If one wanted to deal with the abundancy index directly, using the Hardy-Ramanujan circle method, then one would need to control the integrand
in
$$
\frac{\sigma(n)}{n}\, =\, 
 e^{n y_n+ \ln(F(y_n))}
\int_{-\pi}^{\pi} e^{-in\theta} \exp\left(\ln\left(\frac{F(y_n-i\theta)}{F(y_n)}\right)\right)\, \frac{d\theta}{2\pi}\, ,
$$
where $\kappa_n=1/n$ and $y_n$ solves $-F'(y_n)/F(y_n) = 1/\kappa_n=n$. (So $y_n \sim n$ as $n \to \infty$.) 

There appear to be spikes at the rational multiples of $\pi$ in Figure \ref{fig:Farey}, which is a plot for $y>0$ small.
They are presumably arranged as in the sequence of Farey fractions.
If one were to consider $F(y-i\theta)/F(y)$ at $\theta = 2 \pi f$ for Farey fractions $f$, then the first consideration
would be the limits there, as $y \to 0^+$.
For this purpose determining the Cesaro limit of the abundance index along arithmetic sequences is useful.
A good reference is Apostol, Chapter 3 \cite{Apostol}.

It is interesting that Hardy and Ramanujan used such tools in their paper on their partition formula, 
although later researchers have shown that such are not necessary \cite{HardyRamanujan}.
One may take a single major arc to be a small interval near zero: see for example Newman \cite{Newman}.
We also include our own streamlined version of the argument in Appendix \ref{sec:HR} with the goal
of being self-contained.
But for the abundancy index, one would need the full apparatus of Farey fractions and Ford circles, we are sure.

Farey fractions have been related to mathematical physics for example in references \cite{KlebanOzluk,GuerraKnauf,ContucciKlebanKnauf}.

\section{Motivation}

Abdesselam has proposed an interesting conjecture related to log-concavity of a family of important combinatorial sequences.
Suppose that $\ell \in \N$ and $(\pi_1,\dots,\pi_{\ell})$ is an $\ell$-tuple of permutation of $[n]=\{1,\dots,n\}$ satisfying
that $\pi_i \pi_j = \pi_j \pi_i$ for all $i,j \in \{1,\dots,\ell\}$.
In other words, the permutations commute with one another.
Now let the subgroup generated by them be $\langle \pi_1,\dots,\pi_{\ell} \rangle$.
Then this subgroup, acting on $[n]$ may be used to decompose $[n]$ into some number of orbits $K(\pi_1,\dots,\pi_{\ell})$.
Then let $A(\ell,n,k)$ be the number of $\ell$-tuples $(\pi_1,\dots,\pi_{\ell})$ which commute and such that $K(\pi_1,\dots,\pi_{\ell})=k$.

In \cite{Abdesselam}, he conjectured that
$$
\big(A(\ell,n,k)\big)^2\, \geq\, A(\ell,n,k-1) A(\ell,n,k+1)\, ,
$$
for all $k \in \{2,\dots,n-1\}$. This is the log-concavity condition which is important in combinatorics as well as probability
theory.
He noted that it follows easily for $\ell=1$, from properties of the Stirling numbers of the first kind.
In a paper with coauthors \cite{AbdesselamStudents}, he noted that the $\ell=2$ case is a conjecture
already made by Heim and Neuhauser in \cite{HeimNeuhauser}.

Abdesselam had already proved his conjecture in a number of interesting regimes, when we collaborated on a weak central limit theorem
result \cite{AbdesselamStarr}.
The weak central limit theorem result was not sufficient to establish log concavity.
But in an interesting development, Abdesselam did establish that log-concavity result in that regime by a tour-de-force (in our opinion).
He not only proved a local central limit theorem, but he actually obtained all asymptotics down to order-1, so that any errors in his formulas
are strictly smaller than order 1, namely $o(1)$, as $n \to \infty$.
A similar result had been done for combinatorial numbers equivalent to the sum of the coefficients
by Bringmann, Franke and Heim \cite{BFH}.
But that was a tour-de-force itself (in our opinion).

Let us state an example of one of Abdesselam's results from \cite{Abdesselam}.
For $\ell=2$, we have the numbers of interest in this note.
He proved that, for $k_n$ such that $k_n/\sqrt{n} \to x \in (0,\infty)$
$$
\frac{A(2,n,k_n)}{n!}\, \sim\,
\frac{\exp\left(-2k_n \ln\left(\frac{k_n}{\sqrt{n}}\right)+ 2 k_n + k_n \ln(\zeta(2))+\frac{k_n^2 \Big(-\ln(n)+2\ln\left(\frac{k_n}{\sqrt{n}}\right)-2\ln(2\pi)\Big)}{4\zeta(2)n}
\right)}{2\pi n}\, \, ,
$$
as $n \to \infty$.
He also proved that this is sufficient to establish log-concavity in this regime for sufficiently large $n$.

Our main motivation was to establish this in an alternative way.
But we have only succeeded so far for $\ell=2$.
The difference between $\ell=2$ and higher values of $\ell$ is that for $\ell=2$ the relevant generating function, $F(y)$
is related to Dedekind's eta function.
Therefore, by Euler's pentagonal identity and the Poisson summation formula, one may establish the modular symmetry
that we have stated already (which was discovered by Dedekind, himself).
This makes the result such as the one above relatively easy to establish.
As an alternative to demonstrate the main point: Hardy and Ramanujan's asymptotic formula
for the partition numbers is even easier to derive using the modular symmetry.
We suggest that the interested reader see Appendix \ref{sec:HR}, first.

But Abdesselam did not stop at $\ell=2$.
He proved a number of interesting results, which we believe are of interest in number theory itself.
It remains a long-term goal of ours to fully understand the ramifications of his results.
For now, we constrain our attention to $\ell=2$.

\subsection{Re-derivation of the asymptotics for $\ell=2$}

Note that the asymptotic series for $F$ is 
$$
F(y)\, =\, \frac{\zeta(2)}{y} + \frac{\ln(y)}{2} - \frac{\ln(2\pi)}{2} - \frac{y}{24} + O(y^p)\, ,
$$
for all $p$ as $y \to 0^+$, which we obtain by discarding the correction $F(4\pi^2/y)$,
since it is smaller than $y^p$ for every $p$.
Therefore, the asymptotic series for $F'$ near $0$ is 
$$
F'(y)\, =\, -\frac{\zeta(2)}{y^2} + \frac{1}{2y} - \frac{1}{24} + O(y^p)\, ,
$$
for all $p$ as $y \to 0^+$.
Also,
$$
F''(y)\, =\, \frac{2\zeta(2)}{y^3} - \frac{1}{2y^2} + O(y^p)\, ,
$$
for all $p$ as $y \to 0^+$.
Thus, defining $\Sigma^2(y) = (F''(y)/F(y)) - (F'(y)/F(y))^2$, we have
$$
\frac{F'(y)}{F(y)}\, \sim\, -\frac{1}{y}\, ,\qquad
\frac{F''(y)}{F(y)}\, \sim\, \frac{2}{y^2}\ \text{ and }\ 
\Sigma^2(y)\, \sim\, \frac{1}{y^2}\,  ,\ \text{ as $y \to 0^+$.}
$$
Now we will take $k_n = t n^{1/2}$ so that $\kappa_n = k_n/n = t n^{-1/2}$.
We have $y_n \sim \kappa_n = tn^{-1/2}$ as $n \to \infty$, where $y_n$ is the solution of 
$$
-\frac{F(y_n)}{F'(y_n)}\, =\, \kappa_n\, ,
$$
for each $n$.
But we need more precise asymptotics for $y_n$ since the exponential term in Theorem \ref{thm:BR1} is
$$
k_n \ln(F(y_n))+n y_n\, ,
$$
which equals
$$
k_n \ln\left(\frac{\zeta(2)}{y_n}\right) + k_n\ln\left(1+\frac{y_n \ln(y_n)}{2\zeta(2)}-\frac{y_n \ln(2\pi)}{2\zeta(2)}-\frac{y_n^2}{24 \zeta(2)}\right)
+n y_n + O(n^{-p})\, ,
$$
for all $p$ as $n \to \infty$.
We only want order-1 asymptotics. Therefore, we may rewrite this as
$$
k_n \ln\left(\frac{\zeta(2)}{y_n}\right) + k_n\ln\left(1+\frac{y_n \ln(y_n)}{2\zeta(2)}-\frac{y_n \ln(2\pi)}{2\zeta(2)}\right)
+n y_n + o(1)\, ,
$$
as $n \to \infty$.
Let us write $y_n = \kappa_n Y_n = tn^{-1/2} Y_n$, where we note $Y_n \sim 1$ as $n \to \infty$.
Then the quantity of interest is 
$$
k_n \ln\left(\frac{\zeta(2)}{\kappa_n}\right) - k_n \ln(Y_n) + k_n \kappa_n \left(\frac{\ln(\kappa_n)}{2\zeta(2)}-\frac{\ln(2\pi)}{2\zeta(2)}\right)
+n \kappa_n Y_n + o(1)\, ,
$$
as $n \to \infty$.
Therefore, we only need to keep the asymptotics for $Y_n$ to order $n^{-1/2}$.
But also writing $Y_n = 1+\widetilde{Y}_n$ we have
$-k_n\ln(Y_n)+n\kappa_n Y_n\, =\, n\kappa_n + n \kappa_n \widetilde{Y}_n-k_n\ln(1+\widetilde{Y}_n)$
and $n\kappa_n \widetilde{Y}_n-k_n\ln(1+\widetilde{Y}_n) = O(n^{1/2} \widetilde{Y}_n^2)$.
But
$$
\widetilde{Y}_n\, =\, - \frac{\kappa_n}{2\zeta(2)}\, \left(1+ \ln(\kappa_n) - \ln(2\pi)\right) + o(n^{-1/2})\, .
$$
So we see that $n^{1/2} \widetilde{Y}_n^2=o(1)$ as $n \to \infty$.
So the exponential quantity of interest equals
$$
k_n \ln\left(\frac{\zeta(2)}{\kappa_n}\right) + k_n \kappa_n \left(\frac{\ln(\kappa_n)}{2\zeta(2)}-\frac{\ln(2\pi)}{2\zeta(2)}\right)
+n \kappa_n + o(1)\, ,
$$
as $n \to \infty$.
Therefore, we have derived the following consequence of Theorem \ref{thm:BR1}.
\begin{prop}
Suppose that $k_n$ is a sequence such that $k_n/\sqrt{n}$ is bounded.
Then
$$
\frac{k_n! A(2,n,k_n)}{n!}\, \sim\,
\frac{\exp\left(-k_n \ln(\kappa_n) + k_n + k_n \ln(\zeta(2))+\frac{k_n\kappa_n \ln(\kappa_n)}{2\zeta(2)}-\frac{k_n\kappa_n \ln(2\pi)}{2\zeta(2)}
\right)}{\sqrt{2\pi n/\kappa_n}}\, \, ,
$$
as $n \to \infty$.
\end{prop}
This proposition was already proved by Abdesselam as a small part of his much more general set of results in \cite{Abdesselam}.
See his Theorem 1.1.
\begin{remark}
The number of compositions of $n$ with $k$ parts is $\binom{n-1}{k-1}$.
The above formula suggests that, for $k_n \sim x \sqrt{n}$, the average value of $\prod_{j=1}^{k} \left(\frac{\sigma(\nu_j)}{\nu_j}\right)$ among the choices of $(\nu_1,\dots,\nu_k)$, 
is given by the product of two terms, asymptotically.
The first term is $\zeta(2)^k$, as expected (because the Cesaro limit
of the abundance index is $\zeta(2)$).
The second factor is 
$\exp\left(-\frac{x^2}{4\zeta(2)}\,\ln(n)+\frac{x^2}{2\zeta(2)}\,\ln(x)
-\frac{x^2}{2}\left(1+\frac{\ln(2\pi)}{\zeta(2)}\right)\right)$.
\end{remark}
%

\subsection{Proof of Corollary \ref{cor:logConcave}}

We have established the Bahadur-Rao type formula in Theorem \ref{thm:BR1} in a manner which allows us to vary $k_n$.
Replacing $k_n$ by $k_n+1$ and/or $k_n-1$ and Taylor expanding leads to that result.

Note that if
$$
-\frac{F(y_n)}{F'(y_n)}\, =\, \frac{k_n}{n}\, ,
$$
and if we define $\epsilon_{n}^{\pm}$ such that
$$
-\frac{F(y_n+\epsilon_n^{\pm})}{F'(y_n+\epsilon_n^{\pm})}\,  =\, \frac{k_n \pm 1}{n}\, ,
$$
then a Taylor expansion gives
$$
 \left(\frac{F(y_n)F''(y_n)}{\big(F'(y_n)\big)^2}-\frac{F'(y_n)}{F'(y_n)}\right) \epsilon_n^{\pm} + O\Big(\big(\epsilon_n^{\pm}\big)^2\Big)\,
=\, \pm \frac{1}{n}\, .
$$
So this means
$$
\epsilon_n^{\pm}\, \sim\, \pm\frac{1}{n \big(\mathcal{K}(y_n)\big)^2 \mathcal{V}(y_n)}\, .
$$
Now we can decompose
$$
2k_n \ln\big(F(y_n)\big) +2 n y_n - \sum_{\sigma \in \{+,-\}} \Big((k_n+\sigma \cdot 1) \ln\big(F(y_n+\epsilon_n^{\sigma})\big)+n(y_n+\epsilon_n^{\sigma})\Big)\, ,
$$
as
$$
 \mathcal{A}_n + \mathcal{B}_n\, ,
$$
where
$$
\mathcal{A}_n\, =\, - \ln\big(F(y_n+\epsilon_n^+)\big)+\ln\big(F(y_n+\epsilon_n^-)\big)\, ,
$$
and
$$
\mathcal{B}_n\, =\, \sum_{\sigma \in \{+,-\}} \Big(-k_n \ln\big(F(y_n+\epsilon_n^{\sigma})\big)+k_n \ln\big(F(y_n)\big)-n\epsilon_n^{\sigma}\Big)\, .
$$
By Taylor expanding to first order, we have
$$
\mathcal{A}_n\, =\, -\frac{F'(y_n)}{F(y_n)} \left(\epsilon_n^+-\epsilon_n^-\right)\, \sim\, \frac{1}{\mathcal{K}(y_n)}\, \cdot \frac{2}{n \big(\mathcal{K}(y_n)\big)^2\mathcal{V}(y_n)}\, ,
$$
since $-F'(y_n)/F(y_n) = 1/\mathcal{K}(y_n)$.
So
$$
\mathcal{A}_n\, \sim\, \frac{2}{n \big(\mathcal{K}(y_n)\big)^3 \mathcal{V}(y_n)}\, .
$$
So we will complete the proof if we show that $\mathcal{B}_n$ is asymptotic to $-1/\Big(n \big(\mathcal{K}(y_n)\big)^3 \mathcal{V}(y_n)\Big)$.

For $\mathcal{B}_n$ we need to Taylor expand out to second order. For a general $\delta$ small
$$
-k_n \ln\big(F(y_n+\delta)\big)\,
=\, -k_n \ln\big(F(y_n)\big) - k_n\, \frac{F'(y_n)}{F(y_n)}\, \delta - \frac{k_n}{2}\, \mathcal{V}(y_n)\, \delta^2 + O(\delta^3)\, .
$$
But note that $F'(y_n)/F(y_n)=-1/\kappa_n=-n/k_n$.
So
$$
-k_n \ln\big(F(y_n+\delta)\big)\,
=\, -k_n \ln\big(F(y_n)\big) + n \delta - \frac{k_n}{2}\, \mathcal{V}(y_n)\, \delta^2 + O(\delta^3)\, .
$$
Therefore, we have
$$
-k_n \ln\big(F(y_n+\delta)\big)+k_n\ln\big(F(y_n)\big)-n\delta\, =\, -\frac{k_n}{2}\, \mathcal{V}(y_n)\delta^2 + O(\delta^3)\, .
$$
We can rewrite this as 
$$
-k_n \ln\big(F(y_n+\delta)\big)+k_n\ln\big(F(y_n)\big)-n\delta\, =\, -\frac{n \mathcal{K}(y_n)\mathcal{V}(y_n)}{2}\, \delta^2 + O(\delta^3)\, .
$$
So, summing over $\delta = \epsilon_n^+,\epsilon_n^-$, and using $(\epsilon_n^+)^2=(\epsilon_n^-)^2$, we get
$$
\mathcal{B}_n\,
\sim\, -\frac{n \mathcal{K}(y_n) \mathcal{V}(y_n)}{1} \cdot \frac{1}{n^2 \big(\mathcal{K}(y_n)\big)^4 \big(\mathcal{V}(y_n)\big)^2}\, .
$$
This is the desired formula, so we are done.


\section{Relation to other results}

In \cite{Tripathi}, Raghavendra Tripathi considered the case where $n-k$ is order 1.
This covers the most important part of the regions we cannot cover.
We require in the present note that $k \to \infty$ as $n \to \infty$, and also that $k/n < 1 - \epsilon$ for all $n$,
where $\epsilon>0$ is fixed.
Tripathi covered the most interesting part of the regime where $k/n$ converges to 1.

We note that the Okounkov-Nekrasov polynomials are closely related to the d'Arcais polynomials.
Hong and Zhang considered asymptotics of the coefficients of the Okounkov-Nekrasov polynomials for small $k$ in \cite{HongZhang}.
This is also interesting, although not directly related.

The most interesting cases not yet covered would be $\ell=3,4,\dots$, except of course that Abdesselam has 
a complete asymptotic expansion down to the level smaller than $O(1)$ in the typicality regime \cite{AbdesselamSole}.
That of course is the most important regime.
But it would still be interesting to have complementary results such as the large deviation principle.
The starting point would presumably be Abdesselam's article or else a full explanation of any quasimodularity properties
that may exist for $\ell>2$.

\section*{Acknowledgments}
I am grateful to Malek Abdesselam.

\appendix

\section{Warm-up:  Hardy-Ramanujan partition formula}
\label{sec:HR}
As a warm-up, we will rederive the Hardy-Ramanujan asymptotics, much the same way they did, but without using Farey fractions because
they are not needed for this particularly easy application.
This is important to us because our main result is also proved using the circle method, but without needing Farey fractions.

Start with the generating function
$$
(q;q)_{\infty}^{-1}\, =\, \sum_{n=1}^{\infty} p(n) q^n\, ,
$$
for $|q|<1$.
Therefore, by Cauchy's integral formula, taking a contour integral over a circle of radius $e^{-y}$,
$$
p(n)\, =\,e^{ny} \int_{-\pi}^{\pi} e^{-in\theta}\, e^{F(y-i\theta)}\, \frac{d\theta}{2\pi}\, .
$$
Note that by properties we mentioned before, $F$ is strictly decreasing and convex.
It is easy to see that $\lim_{n \to \infty} F'(y)=0$, while $\lim_{y \to 0^+} F'(y)=-\infty$.

Let $y_n$ be the point such that 
$$
F'(y_n)\, =\, -n\, .
$$
So $y_n \sim \sqrt{\zeta(2)/n}$, as $n \to \infty$.
Note that then $F''(y) \sim 2n^{3/2}/\sqrt{\zeta(2)}$, as $n \to \infty$.
Then, heuristically at least, we expect the following lemma.
\begin{lemma}
\label{lem:HR}
For each $n$, let $y_n$ be the unique point such that $F'(y_n)=-n$. Then
$$
p(n)\, \sim\, e^{ny_n + F(y_n)} \int_{-\infty}^{\infty} e^{-n^{3/2}\theta^2/\sqrt{\zeta(2)}}\, \frac{d\theta}{2\pi}\, =\, 
\frac{e^{-n y_n + F(y_n)}}{2^{5/4} 3^{1/4} n^{3/4}}\, ,
$$
as $n \to \infty$.
\end{lemma}
\noindent
\begin{proof}
{\em We reiterate that $F'(y_n)=-n$ implies $y_n \sim \sqrt{\zeta(2)/n}$.}
$$
F'(y)\, =\, - \sum_{n=1}^{\infty} \frac{n e^{-ny}}{1-e^{-ny}}\, .
$$
So by a Riemann sum approximation to a Riemann integral
$$
F'(y)\, \sim\, -\frac{1}{y^2}\, \int_0^{\infty} \frac{x e^{-x}}{1-e^{-x}}\, dx\, =\, - \frac{\zeta(2)}{y^2}\, ,
$$
as $y \to 0^+$.
So if $F'(y_n)=-n$ then $n \sim \zeta(2)/y_n^2$ as $n \to \infty$. That means $y_n \sim n^{-1/2} \sqrt{\zeta(n)}$ as $n \to \infty$.

{\em We can simplify by trigonometry.}
$$
F(y)\, =\, -\sum_{n=1}^{\infty} \ln(1-e^{-ny})\, =\, \sum_{n=1}^{\infty} \sum_{k=1}^{\infty} \frac{1}{k}\, e^{-nky}\, =\, \sum_{m=1}^{\infty} \frac{\sigma(m)}{m}\, e^{-my}\, .
$$
So
$$
F(y-i\theta)\, =\, \sum_{n=1}^{\infty} \frac{\sigma(n)}{n}\, e^{-ny} e^{in\theta}\, ,
$$
and
$$
F(y-i\theta) - F(y)\, =\, \sum_{n=1}^{\infty} \frac{\sigma(n)}{n}\, e^{-ny} (e^{in\theta}-1)\, .
$$
This means
\begin{equation}
\begin{split}
\label{eq:FdiffReal}
\operatorname{Re}(F(y-i\theta))-F(y)\, 
&=\, - \sum_{n=1}^{\infty} \frac{\sigma(n)}{n}\, e^{-ny} (1-\cos(n\theta))\\
&=\, -2\, \sum_{n=1}^{\infty} \frac{\sigma(n)}{n} e^{-ny} \sin^2\left(\frac{n\theta}{2}\right)\, .
\end{split}
\end{equation}
It also means
\begin{equation}
\label{eq:FImag}
\operatorname{Im}(F(y-i\theta))\, 
=\, \sum_{n=1}^{\infty} \frac{\sigma(n)}{n}\, e^{-ny} \sin(n\theta)\, .
\end{equation}

{\em We have asymptotics for the real part when $\theta = y_n^{3/2} \Theta$.}
\begin{equation*}
\begin{split}
\operatorname{Re}(F(y-i\theta))-F(y)\, 
&=\, -2\, \sum_{n=1}^{\infty} \frac{\sigma(n)}{n} e^{-ny} \sin^2\left(\frac{n\theta}{2}\right)\\
&\sim\, -\frac{1}{2}\, \sum_{n=1}^{\infty} n\sigma(n) e^{-ny} \theta^2\, ,
\end{split}
\end{equation*}
as $\theta \to 0$, by equation (\ref{eq:FdiffReal}).
But note that
$$
\sum_{n=1}^{\infty} n \sigma(n) e^{-ny}\, =\, \sum_{n=1}^{\infty} \sum_{k=1}^{\infty} k \mathbf{1}_{k\N}(n) n \sigma(n) e^{-ny}\,
=\, \sum_{k=1}^{\infty} k \sum_{m=1}^{\infty} mk e^{-mky}\, =\, \sum_{k=1}^{\infty} \frac{k^2 e^{-ky}}{(1-e^{-ky})^2}\, .
$$
And using a Riemann sum limit to a Riemann integral, we see that means
\begin{equation*}
\begin{split}
\operatorname{Re}(F(y-i\theta))-F(y)\, 
&\sim\, -\frac{\theta^2}{2y^3}\, \int_0^{\infty} \frac{x^2 e^{-x}}{(1-e^{-x})^2}\, dx\, =\, - \frac{\zeta(2) \theta^2}{y^3}\, ,
\end{split}
\end{equation*}
as $\theta \to 0$ and $y \to 0^+$.
This result follows as long as $\theta$ converges to $0$ and $y \to 0^+$, no matter how fast those two things happen relative to each other.
So rescaling $\theta = y^{3/2} \Theta$, we have
$$
\lim_{ y \to 0^+} \left(\operatorname{Re}(F(y-iy^{3/2}\Theta))-F(y)\right)\, =\, -\zeta(2) \Theta^2\, .
$$

{\em We have  imaginary part asymptotics when $\theta = y_n^{3/2} \Theta$.}
\begin{equation*}
\begin{split}
\operatorname{Im}(F(y-i\theta))-n\theta\, 
&=\, \sum_{n=1}^{\infty} \frac{\sigma(n)}{n} e^{-ny} \left(\sin\left(n\theta\right)-n\theta\right)\\
&\sim\, -\frac{1}{6}\, \sum_{n=1}^{\infty} n^2\sigma(n) e^{-ny} \theta^3\, ,
\end{split}
\end{equation*}
as $\theta \to 0$ because $F'(y)=-\sum_{n=1}^{\infty} \sigma(n) e^{-ny}$ and $\sin(\phi)-\phi \sim -\phi^3/6$ as $\phi \to 0$.
But then
$$
\sum_{n=1}^{\infty} n^2 \sigma(n) e^{-ny}\,
=\, \sum_{n=1}^{\infty} \sum_{k=1}^{\infty} n^2 k \mathbf{1}_{k\N}(n) e^{-ny}\,
=\, \sum_{k=1}^{\infty} k^3 \sum_{m=1}^{\infty} m^2 e^{-mky}\,
=\, \sum_{k=1}^{\infty} \frac{k^3 e^{-ky}(1+e^{-ky})}{(1-e^{-ky})^3}\, .
$$
Using a Riemann sum approximation to a Riemann integral, 
$$ 
\sum_{k=1}^{\infty} \frac{k^3 e^{-ky}(1+e^{-ky})}{(1-e^{-ky})^3}\, \sim\, y^{-4} \int_0^{\infty} \frac{x^3 e^{-x}(1+e^{-x})}{(1-e^{-x})^3}\, dx\, =\, 6 \zeta(2) y^{-4}\, ,
$$
as $y \to 0^+$.
This implies that
\begin{equation*}
\operatorname{Im}(F(y-i\theta))-n\theta\, 
\sim\, -\zeta(2) y^{-4} \theta^{3}\, ,
\end{equation*}
as $\theta \to 0$ and $y \to 0^+$.
But this is true as long as $\theta$ converges to $0$ and $y \to 0^+$, no matter how fast those two things happen relative to each other.
So we may rescale $\theta = y^{3/2} \Theta$, our desired scaling for the real part.
Then we have
\begin{equation*}
\operatorname{Im}(F(y-iy^{3/2}\Theta))-ny^{3/2}\Theta\, 
\sim\, -\zeta(2) y^{1/2} \Theta^{3}\, ,
\end{equation*}
as $y \to 0^+$.

Because of all this we see that we have pointwise convergence of the integrand
\begin{equation}
\label{eq:ptwiseConvergence}
e^{-iny_n^{3/2}\Theta}\, e^{F(y_n-iy_n^{3/2}\Theta)}\, \to\, e^{-\zeta(2)\Theta^2}\, ,\ \text{ as $y \to 0^+$.}
\end{equation}

{\em We have uniform upper bounds for the real part.}
\begin{equation*}
\begin{split}
\operatorname{Re}(F(y-i\theta))-F(y)\, 
&=\, - \sum_{n=1}^{\infty} \frac{\sigma(n)}{n}\, e^{-ny} (1-\cos(n\theta))\\
&=\, -2\, \sum_{n=1}^{\infty} \frac{\sigma(n)}{n} e^{-ny} \sin^2\left(\frac{n\theta}{2}\right)\, .
\end{split}
\end{equation*}
Now, using monotonicity (since $\sin^2(\phi)\geq 0$) and the fact that $\sigma(n)/n \geq 1$, we have
$$
\operatorname{Re}(F(y-i\theta))-F(y)\, \leq\, -2 \sum_{n=1}^{\infty} e^{-ny} \sin^2\left(\frac{n\theta}{2}\right)\, .
$$
But it is easily seen (for example using Wolfram Mathematica) that 
$$
\sum_{n=1}^{\infty} e^{-ny} \sin^2\left(\frac{n\theta}{2}\right)\, =\, \frac{1}{4} \cdot \frac{\sin^2(\theta/2) \cosh(y/2)}{\sinh(y/2) \big(\sinh^2(y/2)+\sin^2(\theta/2)\big)}\, .
$$
So
$$
\operatorname{Re}(F(y-i\theta))-F(y)\, \leq\, 
-\frac{\sin^2(\theta/2) \cosh(y/2)}{2\sinh(y/2) \big(\sinh^2(y/2)+\sin^2(\theta/2)\big)}\, .
$$
Now we note that $\sin^2(\theta/2) \leq \theta^2/\pi^2$ for $-\pi<\theta<\pi$ (for example by convexity of $\sin(\theta/2)$ for $0<\theta<\pi$).
Also, $\tanh(y)\leq y$ for $y\geq 0$.
Therefore,
$$
\operatorname{Re}(F(y-i\theta))-F(y)\, \leq\, -\frac{1}{y}\, \left(1-\frac{\pi^2 \sin^2(y/2)}{\pi^2 \sinh^2(y/2)+ \theta^2}\right)\, =\, -\frac{y^{-1} \theta^2}{\pi^2 \sinh^2(y/2) + \theta^2}\, .
$$
Now let us scale $\theta = y^{3/2} \Theta$. Then we have
$$
\operatorname{Re}(F(y-i\theta))-F(y)\, \leq\,  -\frac{\Theta^2}{\pi^2 y^{-2}\sinh^2(y/2) + y\Theta^2}\, .
$$
For sufficiently large $n$ we will have $y_n$ sufficiently small that we may assume $y^{-2} \sinh^2(y/2) \leq 1$.
So then
$$
\operatorname{Re}(F(y-i\theta))-F(y)\, \leq\,  -\frac{\Theta^2}{\pi^2 + y\Theta^2}\, .
$$
Since $|\Theta|<\pi/y^{3/2}$, we have $y<\pi^{2/3}/\Theta^{2/3}$.
(So $1/(a+by)>1/(a+b\pi^{2/3}/\Theta^{2/3})$ for $a,b>0$ and $-1/(a+by)<-1/(a+b\pi^{2/3}/\Theta^{2/3})$.)
We thus obtain
$$
\operatorname{Re}(F(y-i\theta))-F(y)\, \leq\,  -\frac{\Theta^2}{\pi^2 + \pi^{2/3} \Theta^{4/3}}\, .
$$
Thus
$$
\exp\left(\operatorname{Re}(F(y-i\theta))-F(y)\right)\, \leq\, \exp\left( -\frac{\Theta^2}{\pi^2 + \pi^{2/3} \Theta^{4/3}}\right)\, .
$$
This is integrable. Moreover, this is the modulus of the integrand, after rescaling.

{\em So the result follows by DCT.}
\end{proof}

Now, using the lemma, let us finish the proof of the Hardy-Ramanujan asymptotics.
We already noted, earlier in the article, that 
$$
F(y)\, =\, \frac{\zeta(2)}{y} + \frac{\ln(y)}{2} - \frac{\ln(2\pi)}{2} - \frac{y}{24} + O(y^p)\, ,
$$
for all $p$ as $y \to 0^+$. This is by the modular symmetry, obtained by discarding the correction $F(4\pi^2/y)$,
since it is smaller than $y^p$ for every $p$.
Therefore, the asymptotic series for $F'$ near $0$ is 
$$
F'(y)\, =\, -\frac{\zeta(2)}{y^2} + \frac{1}{2y} - \frac{1}{24} + O(y^p)\, ,
$$
for all $p$ as $y \to 0^+$.
Also,
$$
F''(y_n)\, =\, \frac{2\zeta(2)}{y^3} - \frac{1}{2y^2} + O(y^p)\, ,
$$
for all $p$ as $y \to 0^+$.

Since
$$
F'(y_n)\, =\, -\frac{\zeta(2)}{y_n^2}\, \left(1 - \frac{y_n}{2\zeta(2)}+\frac{y_n^2}{24 \zeta(2)} + O(y_n^p)\right)\, =\, -n\, ,
$$
that means
$$
n y_n\, =\, \sqrt{n \zeta(2)} - \frac{1}{4} + o(1)\, ,\ \text{ as $n \to \infty$.}
$$ 
For example, dropping the third term $n y_n^2 = \zeta(2) - \frac{1}{2}\, y_n + \varepsilon$
where $\varepsilon=O(1/n)$. This gives $4ny_n = \sqrt{16n \zeta(2)+O(1)}-1 = 4 \sqrt{n\zeta(2)}-1+O(1/\sqrt{n})$. 

Then that means
$$
F(y_n)\, =\, \frac{\zeta(2)}{\sqrt{\zeta(2)/n}} \left(1-\frac{1}{4 \sqrt{n \zeta(2)}} + O\left(\frac{1}{n}\right)\right)^{-1} + \frac{1}{2}\, \ln\left(\frac{\sqrt{\zeta(2)}}{\sqrt{n}}\right) 
-\frac{\ln(2\pi)}{2} + O\left(\frac{1}{\sqrt{n}}\right)\, .
$$
In other words,
$$
F(y_n)\, =\, \sqrt{n \zeta(2)} + \frac{1}{4} + \frac{1}{4}\, \ln(\zeta(2)) - \frac{1}{4}\, \ln(n) - \frac{1}{2}\, \ln(2\pi) + o(1)\, .
$$
Therefore,
$$
ny_n + F(y_n)\, =\, 2\sqrt{n\zeta(2)}-\frac{1}{4}\, \ln(n) + \frac{1}{4}\, \ln(2^3\cdot 3) + o(1)\, .
$$
Therefore,
$$
e^{ny_n+F(y_n)}\, =\, \frac{1}{2^{3/4} 3^{1/4} n^{1/4}} e^{2\sqrt{n \zeta(2)}} \left(1+o(1)\right)\, .
$$
So by Lemma \ref{lem:HR}, we obtain
$$
p(n)\, \sim\, \frac{1}{2^2 3^{1/2} n}\, e^{2\sqrt{n\zeta(2)}}\, ,\ \text{ as $n \to \infty$.}
$$
That is Hardy and Ramanujan's formula.


\baselineskip=12pt
\bibliographystyle{plain}

\begin{thebibliography}{10}

\bibitem{Abdesselam}
Abdelmalek Abdesselam.
\newblock Log-Concavity with Respect to the Number of Orbits for Infinite Tuples of Commuting Permutations.
\newblock {\em Ann.~Combin.} (2024)

\bibitem{AbdesselamSole}
Abdelmalek Abdesselam.
\newblock Proof of a conjecture by Starr and log-concavity for random commuting permutations.
\newblock {\em Preprint} (2025)
\newblock {\url{https://arxiv.org/abs/2506.06894}}


\bibitem{AbdesselamStudents}
Abdemalek Abdesselam, Pedro Brunialiti, Tristan Doan and Philip Velie.
\newblock A bijection for tuples of commuting permutations and a log-concavity conjecture.
\newblock {\em Res.~Number Theory} (224) 10:45.

\bibitem{AbdesselamStarr}
Abdelmalek Abdesselam and Shannon Starr.
\newblock A central limit theorem for a generalization of the Ewens measure to random tuples of commuting permutations.
\newblock {\em Preprint} (2025).
\newblock {\url{https://arxiv.org/abs/2505.11469}}

\bibitem{Apostol}
\newblock Tom M.~Apostol.
\newblock {\em Introduction to Analytic Number Theory.}
\newblock Springer-Verlag, Berlin 1974.



\bibitem{BanerjeeWilkerson}
Subho Banerjee and Blake Wilkerson.
\newblock Asymptotic expansions of Lambert series and related $q$-series.
\newblock {\em Intern.~J.~Number Theor.} {\bf 13}, no.~8, pp.~2097--2113 (2017).


\bibitem{BFH}
Kathrin Bringmann, Johann Franke, Bernhard Heim.
\newblock
Asymptotics of commuting 
$\ell$-tuples in symmetric groups and log-concavity.
\newblock {\em Res.~number theory} {\bf 10}, 83 (2024).

\bibitem{BryanFulman}
Jim Bryan and Jason Fulman.
\newblock Orbifold Euler characteristics and the number of commuting n-tuples in symmetric groups.
\newblock {\em Ann.~Combin.} {\bf 2}, 1--6 (1998).



\bibitem{ContucciKlebanKnauf}
Pierluigi Contucci, Peter Kleban and Andreas Knauf.
\newblock A Fully Magnetizing Phase Transition.
\newblock {\em J.~Statist.~Phys.} {\bf 97}, 523--539 (1999).

\bibitem{Davenport}
H.~Davenport.
\newblock \"{U}ber numeri abundantes.
\newblock {\em S.~Ber.~Preu\ss.~Akad.~Wiss., math.-nat.~Kl} 830--837 (1933).

\bibitem{DemboZeitouni}
Amir Dembo and Ofer Zeitouni.
\newblock {\em Large Deviations Techniques and Applications. Second Edition.}
\newblock Springer-Verlag, New York (1998).


\bibitem{ErdosPacificJ}
P.~Er\H{o}s.
\newblock On the distribution of numbers of the form $\sigma(n)/n$ and on some related questions.
\newblock {\em Pac.~J.~Math.} {\bf 52}, no.~1, 59--65 (1974).

\bibitem{FlajoletSedgewick}
Phillipe Flajolet and Robert Sedgewick.
\newblock {\em Analytic Combinatorics.}
Cambridge University Press, Cambridge, UK 2009.

\bibitem{GuerraKnauf}
Francesco Guerra and Andreas Knauf.
\newblock Free Energy and Correlations of the Number-Theoretical Spin Chain.
\newblock {\em J.~Math.~Phys.} {\bf 39}, 3188--3202 (1998).

\bibitem{HardyRamanujan}
G.~H.~Hardy and S.~Ramanujan.
\newblock Asymptotic formulae in combinatory analysis.
\newblock {\em Proc.~London Math.~Soc.} {\bf 2}, no.~17, 75--115 (1918).

\bibitem{HeimNeuhauser}
B.~Heim and M.~Neuhauser.
\newblock Horizontal and vertical log-concavity.
\newblock {\em Res.~Number~Theory} {\bf 7}, no.~1, Paper No.~18, 12 pages (2021).

\bibitem{HongZhang}
Letong Hong and Shengtong Zhang.
\newblock Towards Heim and Neuhauser's unimodality conjecture on the Nekrasov-Okounkov
polynomials.
\newblock {\em Res.~Number Theor.} {\bf 7}, article 17, 11 pages (2021).

\bibitem{KlebanOzluk}
P.~Kleban and A.~E.~\"{O}zl\"{u}k.
\newblock A Farey Fraction Spin Chain.
\newblock {\em Commun.~Math.~Phys.} {\bf 203}, 635--647 (1999).

\bibitem{LiStarr}
Timothy Li and Shannon Starr.
\newblock Multifold Convolutions, Generating Functions and 1d Random Walks.
\newblock Preprint 2024.
\newblock{\url{https://arxiv.org/abs/2410.22486}}

\bibitem{Moak}
Daniel S.~Moak.
\newblock The $Q$-analogue of Stirling's Formula.
\newblock {\em Rocky Mountain J.~Math.} {\bf 14}, no.~2 (1984).

\bibitem{MorrillPlatt}
Thomas Morrill and David Platt.
\newblock Robin's inequality for 20-free integers.
\newblock {\em INTEGERS: Electr~J~Combin.~Num.~Theor.} {\bf 21} A28, 7pages (2021).

\bibitem{Newman}
D.~J.~Newman
\newblock A Simplified Proof of the Partition Formula.
\newblock {\em Michigan Math.~J.} {\bf 9}, no.~3, 283--287 (1962).

\bibitem{NiccoloMarkusBenedikt}
Bosio Niccolo\`o, Markus Kuba and Benedikt Stufler.
\newblock Gibbs partitions and lattice paths.
\newblock {\em Preprint}, 2004.
\newblock {\url{https://arxiv.org/abs/2411.03930}}

\bibitem{Pitman}
Jim Pitman.
\newblock {\em Combinatorial Stochastic Processes, Ecole d'Et\'e de Probabilit\'es de Saint-Flour
XXXII-2002.}
\newblock Springer-Verlag, Berlin 2006.


\bibitem{Robin}
Guy Robin.
\newblock Grandes valeurs de la fonctoion somme des diviseurs et hypoth\`ese de Riemann.
\newblock {\em J.~de~Math.~Pures~et~Appl.} Neuvi`eme S\'erie, {\bf 63}, no.~2,187--213 (1987).

\bibitem{Starr}
S.~Starr.
\newblock Some observations about the ``generalized abundancy index.''
\newblock {\em Preprint}, 2025.
\newblock {\url{https://arxiv.org/abs/2505.07051}}

\bibitem{SteinShakarchi}
Elias M.~Stein and Rami Shakarchi.
\newblock {\em Fourier Analysis: An Introduction.}
\newblock Princeton University Press, Princeton, NJ 2003.

\bibitem{StewartTubbenhauer}
Willow Stewart and Daniel Tubbenhauer.
\newblock Representation gaps of rigid planar diagram monoids.
\newblock {\em Preprint}, 2025.
\newblock {\url{https://arxiv.org/abs/2505.05846}}


\bibitem{Tripathi}
Raghavendra Tripathi.
\newblock On log-concavity of the number of orbits in commuting tuples of permutations.
\newblock {\em Res.~Number Theor.} {\bf 10}, article number 78, 9 pages (2024).


\bibitem{Wilf}
Herbert Wilf.
\newblock {\em generatingfunctionology}
\newblock Academic Press, Cambridge, MA 1990.

\end{thebibliography}

\end{document}